# Optimizing wind farms layouts for maximum energy production using probabilistic inference: Benchmarking reveals superior computational efficiency and scalability

**Authors:** Aditya Dhoot[1], Enrico G. A. Antonini[1,2], David A. Romero[1], Cristina H. Amon[1]

**Affiliations:**

1 - University of Toronto, Department of Mechanical and Industrial Engineering, Toronto, ON M5S 3G8, Canada.

2 - Carnegie Institution for Science, Department of Global Ecology, Stanford, CA 94305, USA

**Corresponding author:** Enrico G. A. Antonini, enrico.antonini@mail.utoronto.ca

**Abstract**

Successful development of wind farms relies on the optimal siting of wind turbines to maximize the power capacity under stochastic wind conditions and wake losses caused by neighboring turbines. This paper presents a novel method to quickly generate approximate optimal layouts to support infrastructure design decisions. We model the quadratic integer formulation of the discretized layout design problem with an undirected graph that succinctly captures the spatial dependencies of the design parameters caused by wake interactions. On the undirected graph, we apply probabilistic inference using sequential tree-reweighted message passing to approximate turbine siting. We assess the effectiveness of our method by benchmarking against a state-of-the-art branch and cut algorithm under varying wind regime complexities and wind farm discretization resolutions. For low resolutions, probabilistic inference can produce optimal or nearly optimal turbine layouts that are within 3% of the power capacity of the optimal layouts achieved by state-of-the-art formulations, at a fraction of the computational cost. As the discretization resolution (and thus the problem size) increases, probabilistic inference produces optimal layouts with up to 9% more power capacity than the best state-of-the-art solutions at a much lower computational cost.

# 1 Introduction

In the effort to limit greenhouse gas emissions, curtail the threat of climate change, and achieve energy security, governments are investing heavily in wind farms to harness wind energy. An important steppingstone in facilitating this process is to build fast computational tools that help design efficient large-scale wind farms that maximize power generation and minimize infrastructure costs, while adhering to local land-use, environmental, and mechanical constraints. The wind farm layout optimization (WFLO) is the problem that consists of determining the optimal location of wind turbines within a fixed geographical area to maximize the total power capacity of the wind farm under stochastic wind conditions and wake effects between the turbines. Wake effects created by turbines reduce the wind speed directly downstream of the placed turbines, which decrease the expected power capacity of any turbines placed in their wakes. Additionally, wakes can overlap with each other to further decrease the effective wind speeds.

Existing work on WFLO focuses primarily on two key areas: constructing accurate, computationally efficient wake models and improving layout optimization algorithms. Different approaches exist to model wind turbine wakes, namely analytical and numerical models [1]. Analytical wake models, such as the ones developed by Jensen [2], Larsen [3] and Frandsen *et al.* [4], quickly calculate the wake losses through momentum mixing rate simplifications while neglecting complex flow phenomena occurring in wind farms [5, 6]. On the other hand, computationally expensive models capturing detailed wake interactions can be generated using Reynolds-averaged Navier-Stokes models, large eddy simulations, and vortex-wake models [7, 8, 9, 10, 11, 12, 13]. With regards to the layout optimization, the two most important factors of any wind farm design objective function are minimizing wake effects and maximizing power capacity [14, 15]. A more comprehensive objective function can be based on the wind farm's net-present-value profitability by incorporating installation (i.e., civil and electrical), maintenance, and operations costs in addition to wake effects and power capacity, as illustrated in Ref.

[16, 17]. Furthermore, design constraints such as land-use availability, noise generation, turbine proximity, and other infrastructure constraints generally accompany objective functions [14, 18, 19, 20]. As wind farms increase in size with more turbines [21], the computational tractability of the layout optimization problem has also become a matter of study with the development of new algorithms aimed at reducing the problem's computational cost [22, 23].

Early research on layout optimization algorithms used evolutionary algorithms (EA) to maximize power capacity and minimize infrastructure costs in small wind farms. Mosetti *et al.* [24] formulated the WFLO problem as a discrete-variable problem, modelling the wake interactions using the Jensen's wake model and using a genetic algorithm (GA) to find optimal layouts. With the same problem formulation and by fine-tuning GA parameters, Grady *et al.* [25] obtained better layouts with higher power capacity and lower installation costs. Huang *et al.* [26] further enhanced earlier approaches by applying a distributed genetic algorithm to decompose a large wind farm terrain search space into local search spaces to maximize power capacity. Other applications of EAs can be found in Ref. [19, 27, 28]. Overall, EAs were shown to be effective in finding solutions in this non-linear, non-convex, and non-differentiable optimization problem. Other optimization metaheuristics, such as particle swarm optimization (PSO), ant colony optimization (ACO), among others, have also been applied to this problem. For instance, Chowdhury *et al.* [29, 30] applied PSO to varying WFLO objective functions to handle discrete and continuous variables with non-convex objective functions. Similarly, Wan *et al.* [31] applied a Gaussian PSO while also incorporating a local search strategy based on differential evolution to enhance the optimization results. Hou *et al.* [32, 33, 34] used a PSO to minimize the levelized production costs of a wind farm both in restricted areas and by considering cabling cost.

By discretizing the wind farm terrain into grid cells and representing each cell as a mathematical decision variable, mathematical programming models can be developed for the WFLO problem [35]. The relationship between the decision variables can be encoded as an interaction matrix within the objective function, and hard constraints can be developed for budget, spacing, terrain usage, noise, and other infrastructure or environmental constraints. Integer problems can be solved exactly using well-designed algorithms [36] that are generally implemented in off-the-shelf commercial solvers, such as CPLEX (IBM Corp., Armonk, NY). However, in many cases, models need to be simplified to make the problem linear, tractable and convex [37]. For instance, Donovan [38, 39] and Archer *et al.* [40] used a mixed integer programming (MIP) model to minimize the wake interactions between turbines by varying branching strategies within the branch and bound algorithm to reduce optimization time. Turner *et al.* [41] formulated a binary quadratic integer program (QIP) showing that their solutions outperforms EA for complex wind regimes; nevertheless, it can take up to several hours to converge even for the small problem instances used in the literature. These timelines are incompatible with best wind farm design practices, which call for more interactive means of layout design in the context of multidisciplinary design optimization. Hence, these studies underscore the need for algorithms and implementations that can quickly produce good, nearly optimal solutions of the WFLO problem with lower computational runtimes.

More recent advancements coupled computational fluid dynamics (CFD) wake models directly with optimization algorithms to solve the WFLO problem. Kuo *et al.* [42] proposed an algorithm that couples CFD with MIP to optimize layouts on complex terrains. King *et al.* [43, 44] and Antonini *et al.* [45, 46, 47] instead developed gradient-based optimization approaches which used adjoint methods for the gradient calculation. The use of gradient-based methods was driven by the fact that these methods can outperform genetic algorithms in terms of solution quality and computational cost, as demonstrated by Guirguis *et al.* in Ref. [48, 49] using analytical wake models. Although these implementations have the potential of significantly improving how the complex flow phenomena and wake effects are resolved, the inherently high computational cost is currently seen as their main limitation. In other words, CFD-based WFLO formulations sit at the other end of the trade-off between accuracy and computational cost and, based on the current state-of-the-art, are feasible only for final layout design stages.

Here, we focus on developing a novel, fast optimization algorithm for WFLO problem using a discrete-variable formulation intended for quick screening of potential turbine layout candidates during early design stages, and for generating nearly optimal initial solutions for CFD-based WFLO formulations. The novelty of our algorithm relies on the fact that we pose the binary QIP formulation from Ref. [41] as an undirected graph known as a Markov random field, where the global constraints and pairwise wake interactions are encoded within the graph's edges. We solve this new formulation by applying a probabilistic inference method (maximum-a-posteriori, MAP) on the resulting Markov random fields' joint likelihood function to determine the optimal turbine placements. We hypothesize that, with our new approach, we can substantially decrease the computational time/cost incurred to obtain optimal or nearly optimal wind farm layouts. To generalize our results and to focus only on the effectiveness and efficiency of our proposed algorithm, we use the Jensen's wake model and both idealized and realistic wind turbine performance characteristics and wind roses, consistently with established literature [50, 51, 52]. For a performance comparison, we adopt standard benchmark from the WFLO literature [24, 25], specifically selecting two cases representative of the range of problem complexities typically encountered by practitioners [41]. To further study the effectiveness, efficiency and scalability of our formulation, we compare the turbine layouts produced using the proposed method with an exact, computationally exhaustive branch and cut algorithm under varying wind regime complexities, wind farm discretization resolutions, and number of turbines.

**2 Methods**

The objective of the WFLO is to optimally site turbines to maximize the energy generation of a proposed wind farm. The annual energy production (AEP) of the farm is calculated according to Eq. 1:

$$AEP = \sum_{p=1}^{P} \sum_{q=1}^{Q} t_{p,q} P_{total,p,q},  \qquad (1)$$

where $P$ and $Q$ are the number of wind speed and direction bins, respectively, $t_{p,q}$ is the total time that the wind farm experiences each wind speed-direction state, and $P_{total,p,q}$ is the total power produced by the turbines for a given wind speed and direction. Note that in this work we concern ourselves only with the maximization of energy production; the reader is referred to Ref. [14, 19, 53] for approaches that consider installation costs, environmental impact, land-use and other optimization objectives and constraints. In this study, we use two methods to estimate the power generated by the turbines. In the first one (Eq. 2a), we consider the power generated by each turbine to be proportional to the cube of the effective local wind speeds at the turbine hub height. We adopt this modelling choice to directly compare our results to those published in previous works, which have overwhelmingly calculated the turbine power output in a similar fashion [50, 51, 52, 54]. In the second method (Eq. 2b), we use the power as given by the wind turbine's manufacturer. Hence, depending on the application of our methodology, we can have the following two equations for estimating the power generated:

$$P_{total,p,q} = \sum_{k=1}^{K} 0.3 u_k^3, \qquad (2a)$$

$$P_{total,p,q} = \sum_{k=1}^{K} P_k(u_k), \qquad (2b)$$

where $u_k$ is the effective wind speed of the $k$-th turbine of a wind farm with a total of $K$ turbines. In the next sections, we illustrate how the layout optimization problem can be effectively posed as a Markov random field and efficiently solved using probabilistic inference.

*2.1 Wind farm layout optimization mathematical problem*

Here we formulate the wind farm modeling and optimization by discretizing the land being surveyed for turbine siting into smaller areas. The wind farm area is divided into $N$ cells, where each cell can only hold a single turbine and is represented with a binary variable $x_i$, with $i \in I = \{1, \ldots, N\}$. Additionally, to avoid structural damage to the turbines due to wake interactions, the proximity between any two turbines is set to be at least more than 5 times their rotor radius. Because wake interactions are set to be minimized, the proximity constraint is expected to be met by a large amount most of the times.

To estimate the wake effect on a turbine at location $j$, $j \in I$, generated by a turbine at location $i$ upstream of $j$ we use the Jensen's wake model [2]. This is a well-established and widely used model that provides sufficiently accurate estimations of wake losses over flat terrain [53]. It provides the effective wind speed as a function of their relative position, the wake-decay constant ($\alpha$), and the upstream turbine thrust coefficient ($C_T$). The effective wind speed at a given location $j$ due to multiple wakes (see Fig. 1) is calculated in Eq. 3 by assuming that the kinetic energy deficit at $j$ is the sum of all the kinetic energy deficits caused by the individual wakes produced by turbines at locations $i$ upstream of $j$:

$$\left(1 - \frac{u_j}{u_o}\right)^2 = \sum_{i=1, i \neq j}^{N} \left(1 - \frac{u_i}{u_o}\right)^2 x_i, \qquad (3)$$

where $1 - \frac{u_i}{u_o}$ is the wake effect generated on a turbine at location $j$ by a turbine at location $i$ upstream of $j$. Thus, the effective wind speed $u_j$ of a downstream turbine $j$ can then be calculated with Eq. 4:

$$u_j = u_o \left[1 - \sqrt{\sum_{i=1, i \neq j}^{N} \left(1 - \frac{u_i}{u_o}\right)^2 x_i}\right]. \qquad (4)$$

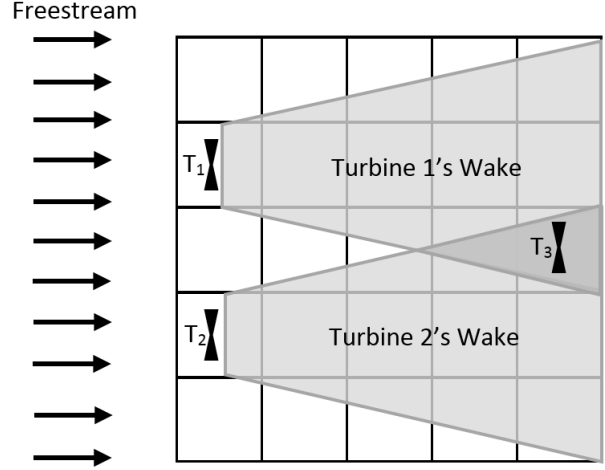

**Figure 1.** Wakes affecting turbine 3 (T3) produced by turbine 1 and 2 (T1, T2).

As illustrated by Turner *et al.* [41], maximizing the wind farm generated power is equivalent to maximizing the effective wind speed at individual turbine locations. In fact, the power generated by a turbine is, to a first approximation, proportional to the cube of the wind speed. In Eq. 5, we then maximize the effective wind speed by only placing turbines at optimal locations $j$ that have an overall smaller kinetic energy deficits and higher effective wind speeds.

$$max \sum_{j=1}^{N} x_j u_o \left[ 1 - \sqrt{\sum_{i=1, i \neq j}^{N} \left(1 - \frac{u_i}{u_o}\right)^2 x_i} \right]. \tag{5}$$

Similarly, a quadratic integer program (QIP) is derived in Ref. [41] by illustrating that maximizing effective wind speed is approximately equivalent to minimizing the kinetic energy losses, which is expressed with Eq. 6:

$$min \sum_{j=1}^{N} \sum_{i=1, i \neq j}^{N} x_j u_o \left(1 - \frac{u_i}{u_o}\right)^2 x_i. \tag{6}$$

An interaction matrix $W \in R^{N \times N}$ can be generated by calculating, prior to the optimization, the kinetic energy deficits at every location $j \in I$ caused by the turbine at locations $i \in I, i \neq j$ considering the wind speed (bin $p$) and direction (bin $q$) probabilities. Each element of this matrix is given by Eq. 7:

$$w_{ij} = \begin{cases} \sum_{p=1}^{P} u_o \sum_{q=1}^{Q} \frac{t_{p,q}}{T}\left(1 - \frac{u_i}{u_o}\right)^2 & , i \neq j \\ 0 & , i = j \end{cases} \quad (7)$$

where $t_{p,q}/T$ represents the probability of occurrence of each combination of wind speed and direction over the entire observation period, $T$. Note that this matrix fully characterizes the complexity of the WFLO problem. Its sparsity pattern is influenced by the wake parameters, the statistical and spatial distribution of the wind resource, and the number of turbines to be placed in the domain, while its size is determined by the discretization of the wind farm domain. All these aspects are typically defined prior to the layout optimization task. For instance, a target wind farm site with a wind regime in which wind can come from any direction with equal probability results in a full, symmetric $W$ matrix, while a site with a unidirectional wind regime results in a non-symmetric, sparse matrix. Practical applications with real-world wind resources will typically result in $W$ matrices that lie between these two extremes. Importantly, the elements of the $W$ matrix, for any given real-world wind resource, and in the context of wind farm design practice, can be calculated with alternative wake models, such as Larsen's [3], Frandsen *et al.*'s [4], or even with CFD simulations of the three-dimensional, turbulent flow in a complex terrain [45, 46, 55]. The $W$ matrix, regardless of how its elements are calculated, is the mathematical expression of the wake- and terrain-induced turbine interactions, and it fully encapsulates the mathematical behavior of the layout optimization problem when the optimization objective is maximizing energy production. In general terms, the QIP formulation can therefore be written with Eq. 8:

$$\text{minimize } X^T W X, \quad (8a)$$

$$\text{subject to } \sum_{j=1}^{N} x_j = K \text{ or } e^T X = K, \quad (8b)$$

$$x_j \in \{0,1\}, \quad (8c)$$

where $X = (x_1, \dots, x_N)$ is a binary decision vector such that $X \in \{0,1\}^N$. Due to budget constraints, land feasibility, government regulations, and grid capacity, the number of turbines to be placed on the wind farm is typically determined prior to the farm development. Thus, a constraint is included in the problem formulation (Eq. 8b), which enforces that only a set number of turbines $K$ will be placed on the farm, where $e$ is a vector $\{1, ..., 1\}^N$. Depending on the behavior of the wake model, the complexity of the wind regime, and the problem state space (i.e., number of turbines and discretization resolution), the resulting QIP can be multi-modal, highly non-convex and non-tractable.

Note that the problem formulation does not require any assumptions regarding the shape of the wind farm area or its subdivision. In Fig. 1, we plotted a squared region divided into squared cells to clearly illustrate the problem's decision variables (the turbine locations). However, the individual cells that compose the wind farm area do not have to form a square or rectangle, they do not have to be numbered in any particular way, they do not even have to be contiguous or form a partition of a space. The shape, dimension, and arrangement of the wind farm area and its cells are chosen depending on the specific wind farm site. As such, problem constraints that are site-specific (e.g., land-use, setbacks) or arise from other disciplines (e.g., noise or environmental requirements) could be incorporated either through appropriate discretization of the wind farm area or as additional terms in the augmented Lagrangian (Eq. 9), to extend our formulation to multi-disciplinary applications.

*2.2 Wind farm layout optimization as a Markov random problem*

Here we pose the binary QIP (Eq. 8) as an undirected graph $G = (V, E)$ known as a Markov random field $V$ is a set of vertices representative of a decision vector, $X \in \{0,1\}^N$ (also known as random variables), connected by edges, $E$, that encode the pairwise wake interactions from $W$ and the turbine number constraint. We start by rewriting Eq. 8 as an augmented Lagrangian in Eq. 9. The Lagrangian encodes the turbine number constraint from Eq. 8b into a quadratic penalty function with penalty factor $\beta$:

$$f(X, \beta) = X^T W X + \beta (e^T X - K)^T (e^T X - K), \qquad (9)$$

which can be further developed into Eq. 10:

$$f(X,\beta) = X^TWX + \beta X^T \begin{bmatrix} 1 & 2 & \cdots & 2 \\ 0 & 1 & \ddots & \vdots \\ \vdots & \vdots & \ddots & 2 \\ 0 & 0 & \cdots & 1 \end{bmatrix} X + \beta(-2K+1)e^TX + K^2. \tag{10}$$

If we omit the constant term $K^2$ as it does not affect the optimization, Eq. 11 shows the unconstrained binary QIP where the parameter $\beta$ controls the smoothness of the penalty contour:

$$\underset{X \in \{0,1\}^N, \beta}{argmin}\ X^T(W + \beta E)X + (-2K\beta + \beta)e^TX. \tag{11}$$

The relationship between the variables from Eq. 11 can be concisely captured using a probabilistic graphical model, $G$ (known as Markov random field), which is a graphical representation of the joint probability distribution of the variables [56]. Probabilistic inference known as maximum-a-posteriori can then be conducted on the resulting Markov random fields' joint likelihood function to determine the optimal turbine placements.

The overall goal of the proposed formulation is to approximate optimal turbine layouts much faster than solving the mathematical model using exact solvers. The key idea behind our proposed approach is that, even though the optimization problem is deterministic, we can use a stochastic approach to solve it more efficiently. Thus, the turbine locations (binary decision variables) are considered as random, with a joint probability distribution between variables (i.e., between turbine locations) that inversely depends on the magnitude of their physical, deterministic interaction as encoded in the $W$ matrix. Once the layout optimization problem is formulated as such, we can leverage an extensive literature in efficient methods for maximization of joint probability distributions to find the turbine positions that are most likely to minimize their interactions and thus have maximum energy production.

Markov random fields (i.e., the probabilistic model of the turbine layout optimization problem) can be decomposed as a product of factors composed of a set of random variables $X_c$ over the set of maximal cliques, $C$, of the graph, $G$ [56]. This decomposition is such that the set of maximal cliques satisfies $X_c \subseteq X\ \forall\ c \in C$ and $\bigcup_{c \in C} X_c = X$ [56]. Intuitively, a clique is a subset of turbines where each of them has a mutual interaction with all the others. Hence, depending on the wake model and the wind regime, the interactions between the random variables may be localized and form maximal cliques. Given a set of binary random variables, $X \in \{0,1\}^N$, the joint probability of a particular configuration of $X$ to occur can be written with Eq. 12:

$$p(X) = \frac{1}{Z} \prod_C \theta_C(X_C), \qquad (12)$$

where $\theta_C(X_C)$ is called potential function and encodes the problem-specific relationships between the random variables, which in this case are the pairwise wake interactions and turbines number constraint. $Z$ is called partition function and is used to normalize the joint distribution. It is expressed as:

$$Z = \sum_X \prod_C \theta_C(X_C). \qquad (13)$$

The potential function, $\theta_C(X_C)$, is generally expressed as an exponential function according to Eq. 14:

$$\theta_C(X_C) = exp(-E(X_c)), \qquad (14)$$

where $E(X_c)$ is called energy function (not to be confused with the wind farm energy generation). The minus sign in Eq. 14 indicates that a random variable assignment with a higher probability has a lower overall energy. The product of the exponential functions can be decomposed in the sum of two sets of energy functions as given in Eq. 15:

$$E(X, \varphi) = \sum_{s \in V} \varphi_s(X_s) + \sum_{(s,t) \in E} \varphi_{st}(X_s, X_t), \qquad (15)$$

where $\varphi_s(X_s)$ represents all unary potentials over the model vertices, $V$ (encoding the decision vector), and $\varphi_{st}(X_s, X_t)$ represents the pairwise potentials over the edges, $E$ (encoding the pairwise wake interactions and the turbine number constraint). If we recall Eq. 11, we see that the terms $(W + \beta_1 E)$ and $(-2K\beta + \beta)$ can be embedded in the coefficients of the quadratic pairwise and linear unary potentials of Eq. 15, respectively. The energy function, $E(X, \varphi)$, can then be interpreted as the original objective function (Eq. 9) to be minimized.

The joint distribution can then be rewritten by incorporating Eq. 15 into Eq. 12:

$$p(X, \varphi) = \frac{1}{Z} exp\left(-\sum_{s \in V} \varphi_s(X_s) - \sum_{(s,t) \in E} \varphi_{st}(X_s, X_t)\right). \tag{16}$$

Determining the configuration of random variables $X_M$ in a graphical model that yields the maximum joint probability (and, consequently, the minimum of the energy function, $E(X, \varphi)$) is known as maximum-a-posteriori inference problem and is formulated as:

$$X_M = \underset{X}{argmax}\, p(X, \varphi); \tag{17}$$

this is equivalent to minimizing the energy function as expressed in Eq. 18:

$$X_M = \underset{X \in \{0,1\}^N}{argmin}\, E(X, \varphi) = \underset{X \in \{0,1\}^N}{argmin} \left(\sum_{s \in V} \varphi_s(X_s) + \sum_{(s,t) \in E} \varphi_{st}(X_s, X_t)\right). \tag{18}$$

Now that the original optimization problem of Eq. 8, reformulated as Eq. 11, has now been cast as a maximum-a-posteriori problem, Eq. 18. Minimizing Eq. 18 is an NP-hard problem, and researchers have focused on approximate minimization algorithms that solve the following, equivalent linear programming (LP) relaxation. By introducing auxiliary marginal variables over the random variables associated with the vertices $\{\mu_s(X_s)\}_{s \in V}$ such that $\sum_{X_s} \mu_s(X_s) = 1$ and for every edge $\{\mu_{st}(X_s, X_t)\}_{(s,t) \in E}$ within the graphical model such that $\sum_{X_s,X_t} \mu_{st}(X_s, X_t) = 1$ [57], the linear program can be defined with Eq. 19 and 20:

$$\underset{\mu \in M(G)}{argmin}\, \varphi^T \mu = \underset{\mu \in M(G)}{argmin} \left(\sum_{s \in V} \sum_{X_s} \mu_s(X_s) \varphi_s(X_s) + \sum_{(s,t) \in \mathcal{E}} \sum_{X_s, X_t} \mu_{st}(X_s, X_t) \varphi_{st}(X_s, X_t)\right), \tag{19}$$

$$M(G) = \left\{\mu \in \mathbb{R}^d | \exists p(X, \varphi)\ subject\ to\ \begin{array}{l}\mu_s(X_s) = \sum_{X_{V \setminus s}} p(X, \varphi) \\ \mu_{st}(X_s, X_t) = \sum_{X_{V \setminus s,t}} p(X, \varphi)\end{array}\right\}, \tag{20}$$

where $M(G)$ is called binary marginal polytope. Eq. 19 and 20 are equivalent to the original problem formulation of Eq. 8 after maximum-a-posteriori probabilistic inference is conducted on the graphical model that encodes decision variables, wake interactions, and turbine number constraint. Below, we discuss a solution strategy for the optimization problem posed in Eq. 19 and 20.

*2.3 Message-passing algorithm*

We hypothesize that approximate solutions can be generated quickly for large Markov random fields by a class of algorithms known as message-passing algorithms that works by exploiting the decomposable factors within the graphical model. Message-passing algorithms approximate maximum-a-posteriori assignments by iteratively passing beliefs locally along the edges of the graphical model in a distributed, decentralized, and asynchronous manner. Intuitively, we can understand this process as follows. Given a candidate turbine layout, i.e., a vector of decision variables that have large joint probability (i.e., a low energy function, thus low turbine interactions), a message passing algorithm would, in any given iteration, modify local groups of decision variables to increase their joint probability (reduce their interaction) regardless of its other neighboring vertices in the graph, and repeat this process asynchronously until convergence.

A widely used message-passing algorithm is the max-product belief propagation (BP) [58], which was shown to produce good empirical results with problem structures similar to Eq. 16, as shown in the fields of computer vision (e.g. Ref. [59, 60]) and computational biology (e.g. Ref. [61]). Because multi-directional wind regimes create cycles within the Markov random field, an algorithm that efficiently computes maximum-a-posteriori configuration for Markov random fields with cycles is tree-reweighted (TRW) message passing [57]. TRW message passing works by decomposing an arbitrary graphical model with cycles into a convex combination of tree-structured distributions to calculate the optimal upper bound while maximizing the lower bound of the energy objective. In this work, we use the sequential tree-reweighted (TRW-S) message passing algorithm developed by Kolmogorov [62] that guarantees that the lower bound will never decrease and, as such, improves its convergence properties. We apply TRW-S to minimize the objective function in Eq. 19 and 20 over the binary marginal polytope to determine the optimal turbine layout configurations. In the Appendix, we provide a description of the TRW-S algorithm as developed by Kolmogorov in Ref. [62, 63].

TRW-S relies on the assumption that the resulting linear program is tight when the integer constraints are relaxed by forming a local polytope $L(G)$, where $M(G) \subseteq L(G)$ and $\min_{\mu \in L(G)} \varphi^T \mu \leq \min_{\mu \in M(G)} \varphi^T \mu$. The relaxed constraints are hardly tight, which requires us to investigate methodologies that can be used to dynamically generate constraints within the message passing paradigm to make the solution space results tight. Sontag *et al.* [64, 65] demonstrated that tighter relaxations can be acquired by iteratively enforcing edge consistency over a small subset of triplet clusters $c \subseteq C$ to generate the polytope $P$. Clusters are chosen at every iteration that improve the dual LP bound to create a sequence of polytopes $P_0 \subseteq P_1 \subseteq P_2 \ldots \subseteq M(G)$ in which the relaxations are continually tighter and approach the marginal polytope [66]. Thus, by iteratively approximating the MIP solution using TRW-S and generating tighter polytopes, we expect that good turbine placements can be acquired in a relatively shorter amount of time compared to other exact solver techniques.

We hypothesize that by applying message passing algorithms to the relaxed WFLO MIP model, while enforcing key constraints using cluster-based approach, optimal or nearly optimal turbine layouts can be generated faster than with traditional exact approaches. In this paper, we conduct a thorough computational study to measure the effectiveness of TRW-S in comparison with branch and cut algorithm to generate turbine layouts under varying wind regime complexities and problem dimensionalities. The specific test cases chosen here reflect typical problem complexities found in wind farm engineering practice, from small wind farms with unidirectional wind regimes to large wind farms under wind regimes with multiple wind speeds and directions, usually showing a dominant wind direction range. Please note that in the following sections, for compactness of our notation, we refer to the proposed TRW-S approximation of the undirected graph/Markov field representation of the underlying QIP simply as the message passing (MP) algorithm.

Note that other inference algorithms exist and could be used to find low-energy solutions. For example, gradient descent is considered one of the simplest optimization techniques, which minimizes the energy by changing independent subsets of variables to find lower-energy configurations [67]. Gradient descent however tends to get stuck in local minima, and stochastic gradient descent are typically implemented to overcome this limitation. Another method that was widely used for MRF inference problems is simulated annealing [68]. However, these methods have been superseded by the newer graph cuts and loopy belief propagation techniques that have proven to be more powerful [69]. Because TRW-S was shown to give consistently strong results [70], it has been chosen as the preferred method also in this study.

*2.4 Optimization workflow*

The workflow of the developed optimization framework is illustrated in Fig. 2 and can be conceptualized as follows. A set of boundary conditions are given during the early design stages of a wind farm project: these include (i) the wind resource of the site, expressed with a wind rose, (ii) the turbine characteristics and performance curves, (iii) site details such as land area, land availability, and environmentally-related setbacks, and (iv) the number of turbines to be placed. A wake model is then used to generate a wake interaction matrix (Eq. 7) by estimating wake effects on the discretized wind farm terrain with the given set of boundary conditions. For example, analytical wake models that have been extensively validated with wind tunnel and field data [71], for use on complex terrains [72], or even CFD simulations [8, 10] could be used for this purpose. In our study, and despite its limitations, we used the Jensen's wake model with the intent of comparing the results of our formulation with standard benchmark problems from the WFLO literature, which has overwhelmingly relied on Jensen's model. We reiterate, however, that any impact of the choice of wake model is reflected in the $W$ matrix, which fully encapsulates all wake interactions and is considered an input to the proposed optimization methodology. Our developed algorithm solves such a defined WFLO problem by optimally siting the turbines in the given available land and by fulfilling the site constraints. Other state-of-the-art optimization algorithms could be used to find optimal layouts, but we will show that our message passing algorithm can generate optimal or nearly optimal layouts at a substantially lower computational cost.

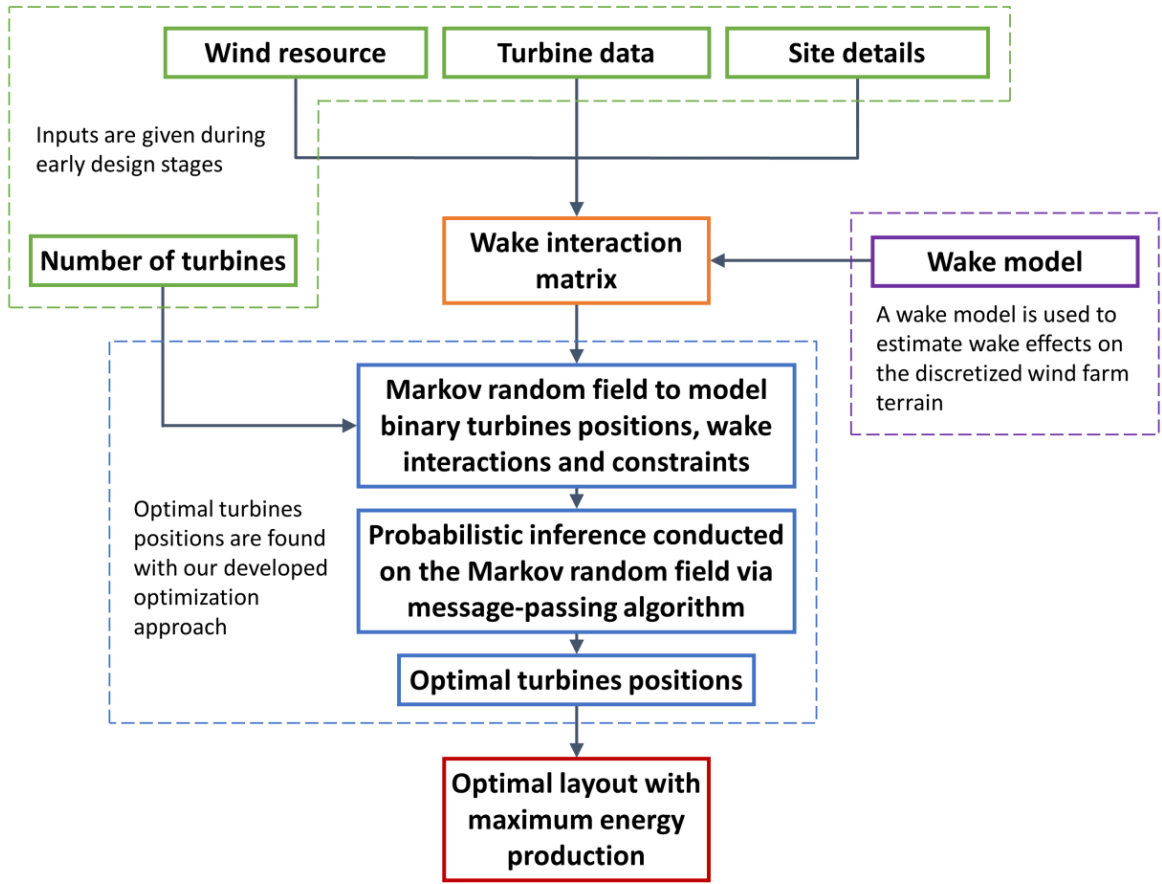

**Figure 2.** Workflow of the developed optimization framework. Our contribution is the development of an efficient and effective optimization approach (in blue) to solve the layout optimization problem and determine an optimal wind farm layout.

*2.5 Numerical experiments*

To study the effectiveness, efficiency and scalability of our formulation, we apply it to several numerical experiments that represent both idealized and realistic wind farm scenarios. We define two benchmark wind resources, WR-1 and WR-36, illustrated in Fig. 3: wind regime WR-1 has a single wind speed and direction state, and WR-36 is a more complex wind regime with 3 wind speeds and 36 directions with varying probabilities. These two wind resource distributions, which result in radically different wake interaction matrices ($W$), allow us to assess the performance of the proposed algorithm in extreme cases, and have been extensively used in the literature to benchmark algorithm performance [24, 25, 41]. In real-world applications, the wind rose for the site would be known, and would be used as an input to the optimization methodology to calculate the corresponding wake interaction matrix.

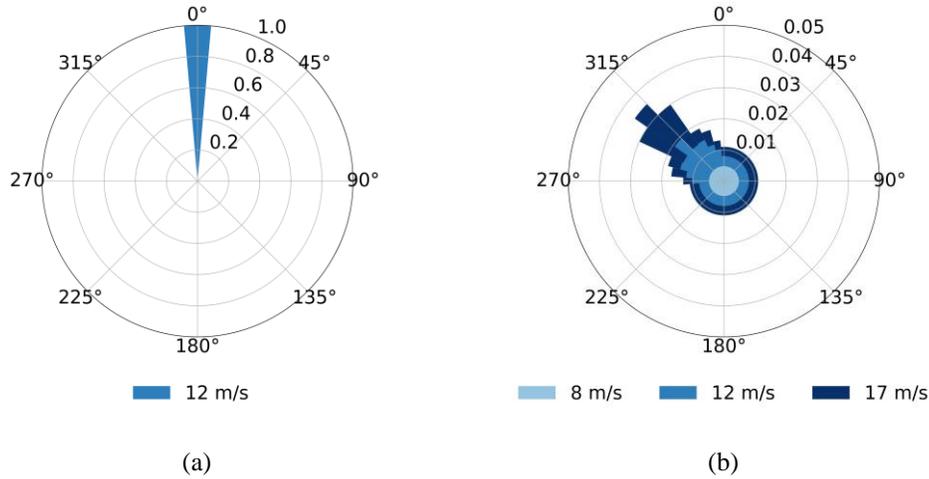

**Figure 3.** Wind roses of unidirectional (a) and multidirectional complex wind regime (b).

Wake modelling, wind farm, and simulation parameters are listed in Tab. 1. When we compare results obtained by our proposed algorithm to those of Ref. [24, 25, 41], we choose their same turbine characteristics, wind farm dimensions, and modelling assumptions to establish a fair comparison. For the more realistic application, we select instead the NREL 5-MW reference wind turbine [73], and we consider a square wind farm area of 49 km$^2$ with varying discrete resolutions (100, 400, and 2,500 square cells). To fully characterize the performance of our algorithm, for each wind rose and discrete resolution, we run a set of cases where we vary the turbine number, $K$.

For the realistic, modern-sized wind farm cases with NREL 5-MW turbines, we compare results obtained by our proposed algorithm with those obtained by the exact state-of-the-art solver CPLEX 12.1 (IBM Corp., Armonk, NY) using a branch and cut algorithm. Initially, the branch and cut algorithm is run using 1 thread with a cut-off of 1 hour ($E_1$); however, if relative optimality is not guaranteed as indicated by its relative gap, further branch and cut iterations are conducted based on 4 threads with a cut-off of 1 hour ($E_2$) and 12 threads with a cut-off of 24 hours ($E_3$). The relative gap serves as a measure of progress toward finding and proving optimality and can be defined as the difference between a best-known solution and a value that bounds the best possible solution. MP simulations are run using source-code in Ref. [74] by applying a total of 5,000 triplet clusters on a single thread with a cut-off of 1 hour. Triplet clusters used to tighten the relaxed marginal polytope (Eq. 20) are generated using modified source-code in Ref. [75] developed from Ref. [66].

**Table 1.** Parameters used in the layout optimization problem.

| Parameters | Values for comparison with results from Ref. [24, 25, 41] | Values for realistic applications |
|---|---|---|
| Farm size | 4 km² (2 km x 2 km) | 49 km² (7 km x 7 km) |
| Turbine rotor radius ($R$) | 20 m | 63 m |
| Turbine hub height ($H$) | 60 m | 90 m |
| Turbine thrust coefficient ($C_T$) | 0.88 | From Ref. [73] |
| Turbine power ($P$) | From Eq. 2a | From Ref. [73] |
| Turbine rated power | n.a. | 5 MW |
| Wake model | Jensen's [2] | |
| Axial induction factor ($a$) | From $C_T$ | |
| Wake-decay constant ($\alpha$) | 0.1 | |
| Exact solver | CPLEX 12.1 (IBM Corp) | |
| Workstation | Intel Xeon Processor - 16 cores @ 2.10 GHz - 256 GB RAM | |

## 3 Results and discussion

In this section, we show the results of our proposed optimization algorithm when applied to the cases described in the previous section, which are useful to understand the effectiveness of our approach. We also compare its results in term of computational time and solution quality to other widely studied algorithms.

*3.1 Comparison with literature results*

In this section, we compare power capacity results from MP with well-established algorithms reported in the WFLO literature, specifically EAs in Ref. [24, 25] and mathematical programming in Ref. [41]. In these cases, the wind farm area is of 4 km$^2$ and divided into 100 cells of 200 m x 200 m. Results in term of power production are reported in Tab. 2 and 3. MP outperforms Mosetti *et al.* [24] for placing 26 and 15 turbines under WR-1 and WR-36, respectively; contrarily, MP performs worse than Grady *et al.* [25] for placing 30 and 39 turbines under WR-1 and WR-36, respectively. Since Mosetti *et al.* and Grady *et al.*'s methods are similar, it follows that Grady *et al.* outperforms the proposed MP thanks to the use of fine-tuned algorithm parameters and a longer run time. Also, results from branch and cut either outperform or match Turner *et al.*'s results [41] even though both methods have similar mathematical formulation.

Figures 4 and 5 show the various layouts for the results presented in Tab. 2 and 3. The algorithms try to find the arrangement that maximizes the power production and, consequently, minimizes the wake losses by moving the turbines further apart. The layouts are a result of different concurring factors, such as the wind speed and direction distributions, the turbine number and performance, and the interspacing constraints. In general, we observe that layouts generated using mathematical programming in Turner *et al.* and CPLEX yield layouts with recognizable spatial structures, while the layouts obtained with EAs in Ref. [24, 25] and MP exhibit less structure but similar energy production values. The results are evidence of the flatness of the objective function landscape near the optimum, since noticeably different layouts result in small differences in power generation values. Importantly, this provides wind farm designers with layout flexibility: multiple layouts may meet energy generation targets, but with markedly different consequences in terms of infrastructure, cost, environmental impact, and other concerns not addressed in this optimization formulation.

**Table 2.** Results comparison under WR-1 with 100 cells.

| K | Model | Power (kW) |
|---|---|---|
| 26 | Mosetti *et al*. [24] | 12,474 |
|  | Turner *et al*. [41] | 12,686 |
|  | CPLEX ($E_1$) | 12,709 |
|  | MP | 12,486 |
| 30 | Grady *et al*. [25] | 14,410 |
|  | Turner *et al*. [41] | 14,410 |
|  | CPLEX ($E_1$) | 14,410 |
|  | MP | 13,972 |

**Table 3.** Results comparison under WR-36 with 100 cells.

| K | Model | Power (kW) |
|---|---|---|
| 15 | Mosetti *et al*. [24] | 13,374 |
|  | Turner *et al*. [41] | 13,671 |
|  | CPLEX ($E_3$) | 13,679 |
|  | MP | 13,395 |
| 39 | Grady *et al*. [25] | 32,377 |
|  | Turner *et al*. [41] | 31,947 |
|  | CPLEX ($E_3$) | 32,818 |
|  | MP | 32,142 |

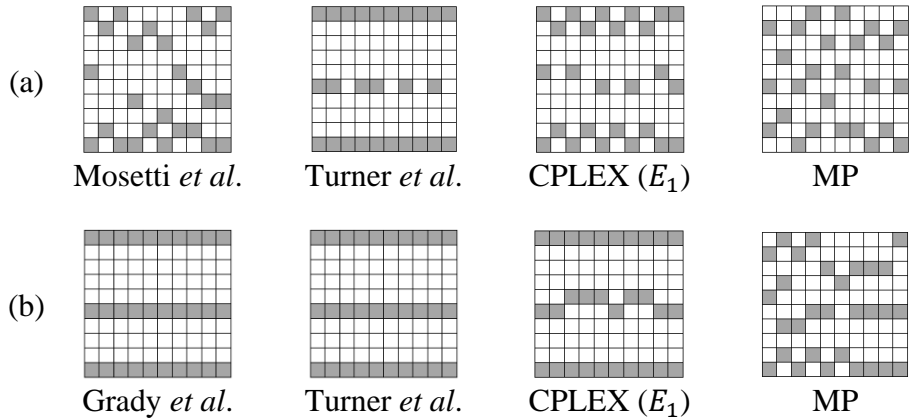

(a) Mosetti *et al*.    Turner *et al*.    CPLEX ($E_1$)    MP

(b) Grady *et al*.    Turner *et al*.    CPLEX ($E_1$)    MP

**Figure 4.** Comparison of layouts for WR-1 scenario with (a) 26 and (b) 30 turbines.

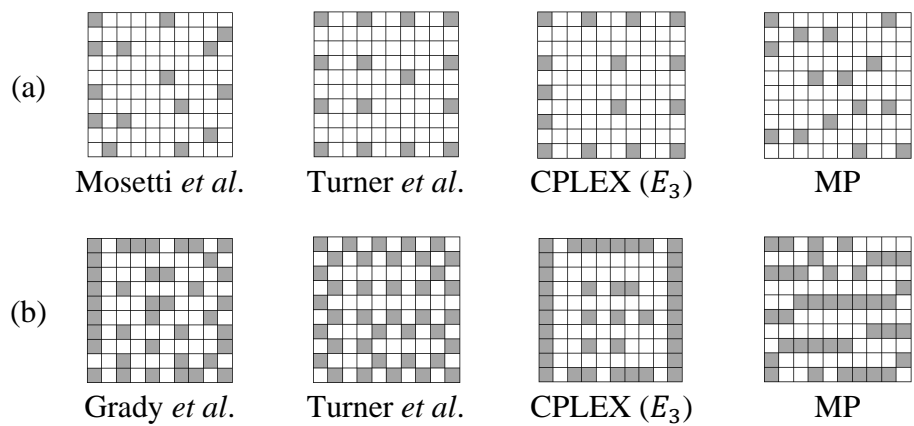

|     |         |         |             |    |
|-----|---------|---------|-------------|----|
| (a) | Mosetti *et al*. | Turner *et al*. | CPLEX ($E_3$) | MP |
| (b) | Grady *et al*.   | Turner *et al*. | CPLEX ($E_3$) | MP |

**Figure 5.** Comparison of layouts for WR-36 scenario with (a) 15 and (b) 39 turbines.

*3.2 Comparison with state-of-the-art solver at 100-cell resolution*

Here, we seek to optimize the layout of a wind farm consisting of NREL 5-MW turbines. The wind farm area is of 49 km² and divided into 100 cells of 700 m x 700 m, which automatically satisfies the turbine proximity constraint. We obtained optimal layouts for varying turbine numbers using branch and cut algorithm and MP methods. Percentage difference in wind farm layout power capacity and ratio of computational times for varying number of turbines for WR-1 and WR-36 are shown in Fig. 6. A positive percentage difference in Fig. 6-a means that MP generated a layout with a higher power capacity. Likewise, a positive time ratio in Fig. 6-b means that MP required a lower computational time than CPLEX. All turbine placements for WR-1 reach optimality using $E_1$ in a few seconds due to the low sparsity of the interaction matrix and the relatively small state space. Under WR-36, the highly dense interaction matrix causes turbine placements to have a larger gap when applying $E_1$, thus, needing application of $E_2$ and $E_3$ to generate better bounds and cuts, requiring up to several days in CPU time to compute optimal solutions. However, a large gap still occurs while placing 20 – 40 turbines due to the relatively large complexity of the turbine constraint polytope as characterized by the difficulty to generate effective cuts and better bounds that well approximate the integer polytope.

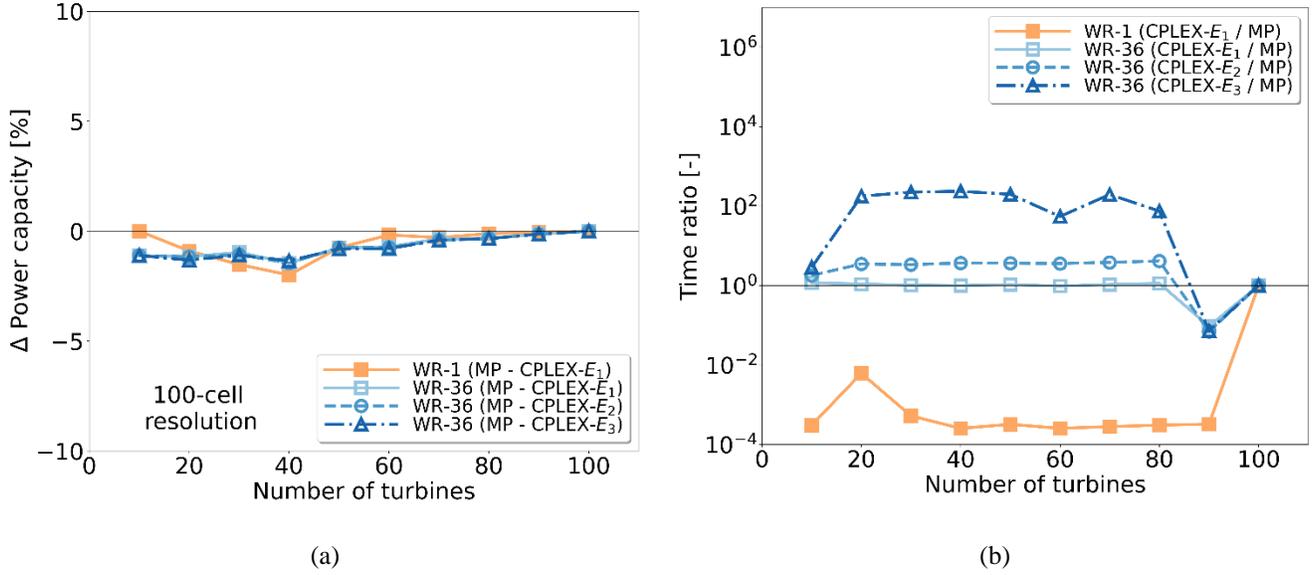

(a)  (b)

**Figure 6.** Percent difference in power capacity (a) and ratio of computational time (b) between MP (our algorithm) and CPLEX for varying turbines under WR-1 and WR-36 with a 100-cell wind farm area. Values above 0% and above $10^0$ indicate that MP produces a layout with, respectively, higher power capacity and lower computational cost.

As shown in Fig. 6-a, many layouts generated using MP capture a lower power capacity compared to $E_1$, $E_2$ and $E_3$ for WR-1 and WR-36, respectively. Under WR-1, MP generates layouts with a power capacity almost equivalent to $E_1$, whereas under WR-36, the minimum percentage difference occurs between $E_1$ and MP, while the maximum occurs between $E_3$ and MP because $E_3$ has a higher cut-off and uses more computing resources to generate similar or better layouts. Power capacity values generated using MP are 0% – 3% lower compared to $E_1$, $E_2$ and $E_3$ for placement of any number of turbines, while varying less than 0.5% between $E_1$, $E_2$ and $E_3$ for any given number of turbines.

MP generates layouts that are consistently within 3% of the power capacity generated using state-of-the-art branch and cut algorithm for a wide variety of cases with varying gaps and cut-off periods. We attribute near optimality to the lack of effective triplet clusters, choice of penalty constant, and local numerical instability caused by locally similar wake decay values. In latter sections, we test again these inefficiencies of MP when generation of cuts and bounds is further challenged when the problem's state space is increased by increasing discrete resolution. In terms of computational time, MP is not capable of efficiently solving the problem given by the domain discretized in 100 cells and under WR-1. However, results in the next section show that, as the wind resource becomes more complex, MP generates nearly optimal layouts in a computational time comparable to $E_1$ or several times faster than $E_2$ and $E_3$.

*3.3 Comparison with state-of-the-art solver at 400-cell resolution*

Resolution is increased by four-fold from 100 cells to 400 cells while keeping the wind farm area to a constant 49 km² and decreasing the square cell size to 350 m x 350 m. Turbine proximity constraint is met without an addition of further constraints due to the adequate distance between centroids of the neighboring cells. Higher resolutions further increase the computational complexity of the problem due to an exponential increase in the state space of the integer program.

We generated optimal layouts for varying turbine numbers using branch and cut algorithm and MP methods for WR-1 and WR-36. As running branch and cut under $E_3$ requires significant computational and memory resources, in some cases, insufficient memory due to large problem size causes the program to terminate without generating a layout. Generation of triplet clusters for 400 cells is a computationally expensive task and can take up to 3 – 4 hours. The graphical model structure based on the marginal polytope generated using the unconstrained objective function (Eq. 19) hardly changes while optimizing for varying number of turbines using MP; thus, we pre-calculate 5,000 triplet clusters based on 10 turbines to save time on layout generation for additional turbines. All results presented using MP are therefore based on pre-generated clusters for individual wind regimes.

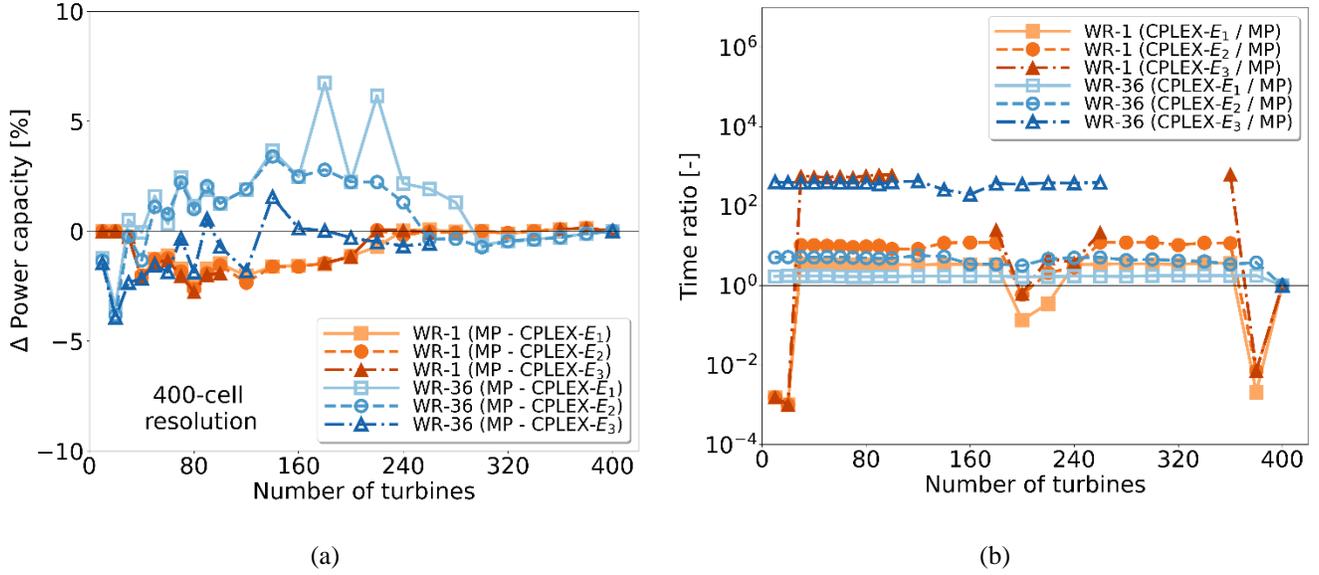

**Figure 7.** Percent difference in power capacity (a) and ratio of computational time (b) between MP (our algorithm) and CPLEX for varying turbines under WR-1 and WR-36 with a 400-cell wind farm area. Values above 0% and above $10^0$ indicate that MP produces a layout with, respectively, higher power capacity and lower computational cost. Missing data points for CPLEX under $E_3$ indicate lack of converged results within the allotted computational resources.

Figure 7 shows the percentage difference in power capacity and ratio of computational times obtained from the three instances of branch and cut and MP for varying number of turbines under both wind regimes. Under WR-1, MP consistently captures layouts with a 0% – 4% lower power capacity than $E_1$, $E_2$ and $E_3$ and with under 1% variation between the three instances. Furthermore, the largest difference occurs in formulations with 30 – 180 turbines due to their difficult turbine constraint polytope. Under WR-36, placing 60 – 280 turbines using MP produces a higher power capacity than $E_1$ and $E_2$. Applying $E_3$ yields tighter bounds and cuts, hence, reducing the gap within this turbine range, which results in $E_3$'s power capacity results that are better than MP. The efficiency of MP in solving the layout problem becomes clear with these cases: MP is 45% to 70% faster than $E_1$, between 5 and 10 times faster than $E_2$, and roughly 3 orders of magnitude faster than $E_3$.

*3.4 Comparison with state-of-the-art solver at 2,500-cell resolution*

Integer optimization of the WFLO problem using any method becomes a magnitude more challenging when the resolution in increased twenty-five-fold to 2,500 cells, while still maintaining the wind farm area to a constant 49 km² and decreasing the square cell size to 140 m x 140 m. At this resolution, placing turbines in adjacent and nearby cells violate turbine proximity constraint. An inequality constraint then is introduced such that if the distance between the two cells is less than 5 times the rotor radius then there can only be a placement of 1 turbine between the two cells. The proximity constraint and 140 m distance between adjacent cells limits the layout feasibility to no more than 280 turbines. Generation of effective bounds, cuts, and triplet clusters is challenging under this case due to the demanding memory requirement and computational complexity. Due to the large state-space we notice that we run out of memory while iterating to tighten the relaxed marginal polytope. Therefore, to preserve computational feasibility we decrease the cut-off time of $E_3$ to 4 hours, remove the generation of triplet clusters, and employ a simple rounding scheme to construct a feasible integral solution for MP.

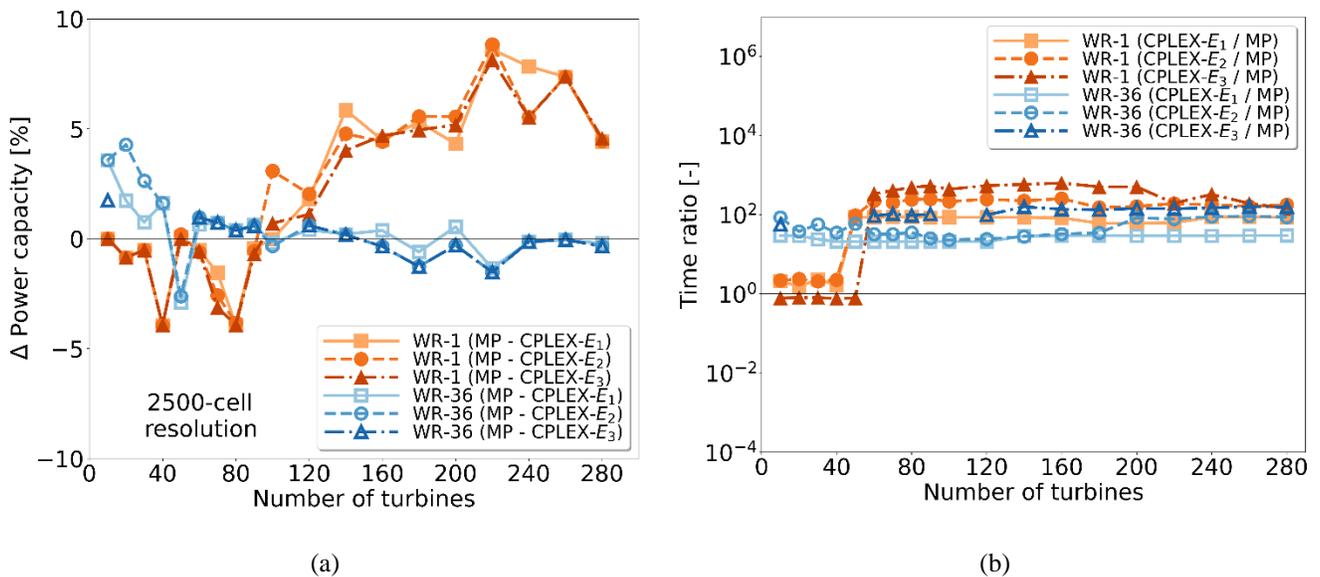

**Figure 8.** Percent difference in power capacity (a) and ratio of computational time (b) between MP (our algorithm) and CPLEX for varying turbines under WR-1 and WR-36 with a 2500-cell wind farm area. Values above 0% and above $10^0$ indicate that MP produces a layout with, respectively, higher power capacity and lower computational cost. Missing data points for CPLEX under $E_3$ indicate lack of converged results within the allotted computational resources.

In Fig. 8 we plot the percentage difference in power capacity and ratio of computational time obtained from the three instances of branch and cut and MP for varying number of turbines under both wind regimes. Results indicate that there is large gap for many placements after applying $E_1$, $E_2$ and $E_3$ for both wind regimes, and $E_3$ generally produces layouts with better power capacity. Even though the lack of triplet cluster generation for MP is expected to produce layouts with uncertain power capacity, we observe that in certain cases power produced by MP is better than $E_1$, $E_2$ and $E_3$. It is conjectured that this occurs due to the difficulty in constructing effective bounds and cuts using CPLEX given the limited amount of resources available and the large problem size. Under WR-1, $E_1$, $E_2$ and $E_3$ produce layouts with up to 5% higher power capacity than MP for up to 100 turbine placements. For higher numbers of turbines, on the other hand, MP produces layouts with up to 9% higher power capacity than branch and cut instances. This occurs as branch and cut is able to generate better bounds and cuts for the turbine constraint polytope with a smaller state space due to a lower demand on computational complexity. Under WR-36, both algorithms produce layouts with similar power capacity for placement of over 50 turbines. Nevertheless, MP performance is approximately better than branch and cut for less than 50 turbines. In terms of computational time, MP generates those solutions to the layout problem several orders of magnitude faster than branch and cut algorithm.

*3.5 Scalability of message-passing algorithm*

In the earlier analysis, we have observed that computational complexity of the integer program increases as the wind farm resolution is increased. While choosing a lower resolution may offer computational feasibility, the resulting turbine configurations may be sub-optimal due to stricter limitations on the placement of the turbines, compared to continuous formulations or discrete formulations with a higher resolution. Therefore, selecting a problem resolution *a priori* while maintaining a handle on tractability poses a challenge for WFLO practitioners.

To illustrate this point, Fig. 9 shows the power capacity for various number of turbines for three cases (100, 400, and 2,500 cells) produced using branch and cut and MP algorithms under WR-1 (left) and WR-36 (right) wind regimes. Under the WR-1 wind regime, both algorithms result in a higher power generation when higher resolutions are used to discretize the wind farm domain. In contrast, problem resolution seems to have no effect on the best achievable power generation values for a given number of turbines under the WR-36 wind regime. If we remember that problem resolution has a direct effect on computational complexity and cost, the complete set of results presented in this work are clear evidence of the complex interaction between the statistical distribution of the wind resource, the problem representation and solution space (discretization resolution), the desired number of turbines, and the computational cost and time required to find a nearly-optimal or optimal solution.

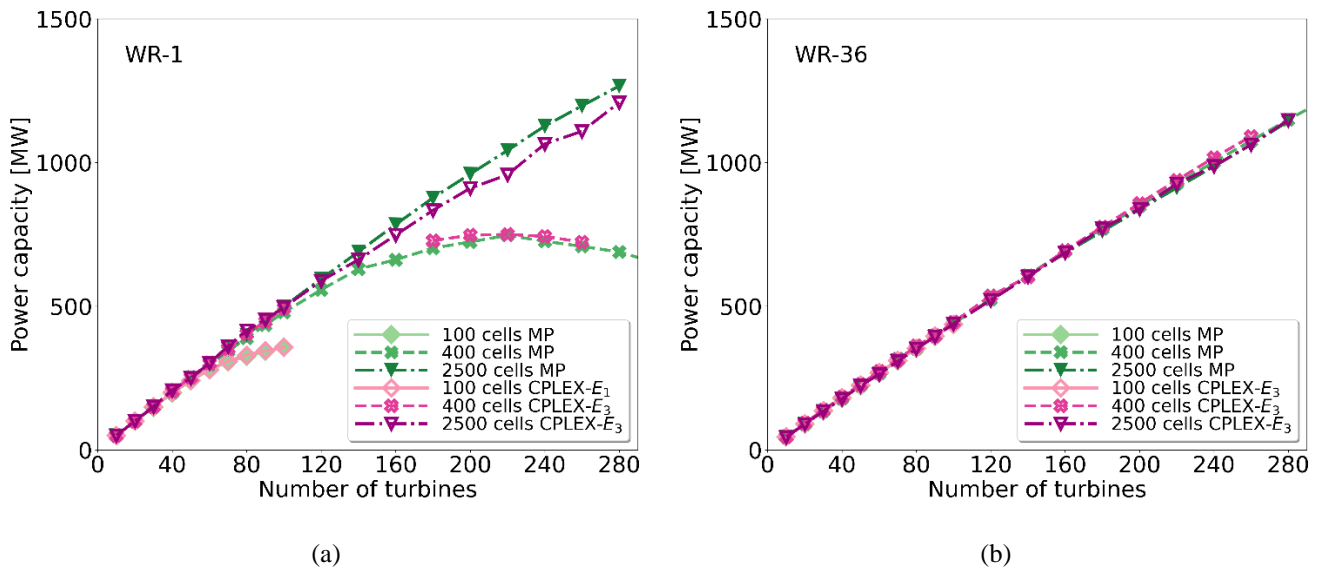

**Figure 9.** Power capacity generated using MP (our algorithm) and CPLEX for 100, 400, and 2,500 square cells under WR-1 (a) and WR-36 (b) wind regimes. CPLEX solutions are taken with $E_3$ settings, except for 100 cells with WR-1, where only $E_1$ is available. Missing data points for CPLEX under $E_3$ indicate lack of converged results within the allotted computational resources.

Scalability of branch and cut and MP poses a challenge due to the difficulty in generating effective cuts and triplets. Nevertheless, for many turbine placements MP's performance is equivalent or better than $E_1$, $E_2$ and $E_3$, particularly at higher resolutions with 400 and 2,500 cells under both studied wind regimes. Interestingly, for these resolutions, both methods produce layouts that are nearly optimal, as evidenced by the large relative gaps between the MIP bounds and the best obtained solution. Hence, in the context of the problem cases studied here, which reflect the complexity of typical WFLO problems found in practice, MP provides a more competitive and scalable approach to generating optimal or nearly optimal layouts in comparison with exhaustive branch and cut while requiring fewer computational resources. Furthermore, nearly optimal layouts generated using MP can be augmented by feeding the layouts to greedy problem-specific heuristics or gradient-based methods to converge to the nearest optimum of the objective landscape. Therefore, MP provides an interesting alternative perspective to generating fast and approximate layouts to the challenging WFLO problem, one more step towards interactive, computer-aided wind farm layout optimization.

## 4. Conclusions

In this paper, we developed a novel, fast optimization algorithm for WFLO problem intended for quick screening of potential turbine layout candidates during early wind farm design stages, or for generating nearly optimal initial solutions for CFD-based WFLO formulations. To facilitate the computational tractability of the problem, generalize the optimization results, and focus on the effectiveness, efficiency and scalability of our algorithm, we discretized the wind farm domain, used the Jensen's wake model, and assumed both idealized and realistic wind turbine performance characteristics and wind roses. We modeled the binary QIP formulation of the WFLO problem as an undirected graphical model by incorporating pairwise wake interactions and turbine number constraints within the graph's edges. Generally, the original QIP formulation requires an exhaustive and computationally expensive branch and cut algorithm to determine optimal layouts. This solver, in many cases, produces sub-optimal layouts with a large relative gap for high discretization resolutions and challenging turbine budget constraints. However, we demonstrated that by using our proposed message passing algorithm on the new undirected graphical model of the WFLO, we can decode good, nearly optimal layouts for a given number of turbines in a substantially lower computational time. To assess the effectiveness, efficiency and scalability of our optimization algorithm, and consistently with established literature, we have conducted a thorough computational study and we have compared it with branch and cut algorithm in generating turbine layouts under varying wind regime complexity and problem resolution relying on commonly used benchmark problems.

The branch and cut algorithm was effective at finding optimal layouts for a low-resolution problem (100 cells) under the two wind regimes tested in this work, while in comparison, MP produced nearly optimal layouts that are consistently within 3% of the power capacity achieved by branch and cut. For test cases with higher resolutions, both methods have difficulty generating cuts to tighten the relaxed polytope and the generated layouts are sub-optimal and rely on locally optimal rounding schemes. At these higher resolutions with 400 and 2,500 cells, MP can outperform branch and cut methods by up to 9%. Additionally, we show that MP can produce better approximate layouts in a shorter period of time. With the exception of the simplest test case tested with 100 cells and unidirectional wind rose (WR-1), MP generated layouts significantly faster, accelerating the optimization process by at least 50%, and it was in many cases several orders of magnitude faster that branch and cut. Moreover, for test cases with higher resolutions and a larger number of turbines, i.e., cases with large number of decision variables, tighter constraints and challenging computational tractability, branch and cut was unable to reach a solution. MP, on the other hand, was able to provide solutions at a lower computational cost.

Overall, MP offers a competitive and scalable alternative to computationally expensive branch and cut algorithms, especially when design engineers are seeking to generate approximate and nearly optimal layouts quickly. For simple wind regimes and low resolutions, branch and cut algorithms still offer a very efficient approach, providing optimal solutions quickly. However, as the wind regimes become more complex or a higher problem resolution is required to better search the solution space, the proposed

MP has a substantial advantage in computational time while providing solutions that are comparable or better than state-of-the-art formulations and, in the most complex cases, providing the only solutions that we could find. The proposed MP can be efficiently used for quick screening of potential turbine layout candidates during early design stages or for generating nearly optimal initial solutions for CFD-based WFLO formulations.

We foresee that a population of approximate layouts generated from MP can be further optimized for other application-relevant objectives (e.g. land-use, noise generation, infrastructure costs) using fine-tuned stochastic methods, evolutionary search methods, or even gradient-based methods. Additionally, our optimization algorithm, which minimizes kinetic energy losses as a proxy for maximizing energy production, can easily accommodate more advanced wake models and more realistic details of a wind farm project, such as wind rose and turbine performance characteristics. To this end, suitable modifications would be needed for the calculation of wake interactions and annual energy production, but the proposed approach would still be valid, since the optimization formulation relies on a wake interaction matrix calculated *a priori* and the annual energy production calculated *a posteriori*.

The use of message passage algorithms to solve the probabilistic inference problem could also enable a more robust wind farm layout optimization with a more comprehensive uncertainty quantification. In fact, recent studies [76, 77] have shown how to compute a measure of uncertainty associated with the graph cut solutions. They showed how the min-marginal energies associated with the label assignments of a random field can be computed using newly developed algorithms based on dynamic graph cuts. With these methods, it is then possible to give confidence bounds on the farm annual energy production and a probability that each determined turbine location is indeed optimal under the wake interaction matrix. This confidence levels on the energy production can be complemented by quantifying uncertainty of both input parameters [78, 79, 80] and wake models [81, 82], which are likely to have a more significant impact in the resulting optimal layout. A comprehensive robust optimization approach that considers and quantifies all these uncertainties, although beyond of the scope of the present study, would be a worthy goal for future work in this area.

**Funding**

This research was funded by the Natural Sciences and Engineering Research Council of Canada (NSERC), grant numbers 397899-2011, 437325-2012 and 462056-2014.

**Appendix: Sequential tree-reweighted (TRW-S) message passing algorithm**

In this work, we use the sequential tree-reweighted (TRW-S) message passing algorithm developed by Kolmogorov [62, 63, 69] to solve the layout optimization problem posed as a Markov random field. We apply TRW-S to minimize the objective function in Eq. 19 and 20 over the binary marginal polytope to determine the optimal turbine layout configurations. In Fig. 10, we present a flowchart of the TRW-S message passing algorithm by illustrating its fundamental steps for finding the optimal solution to the problem.

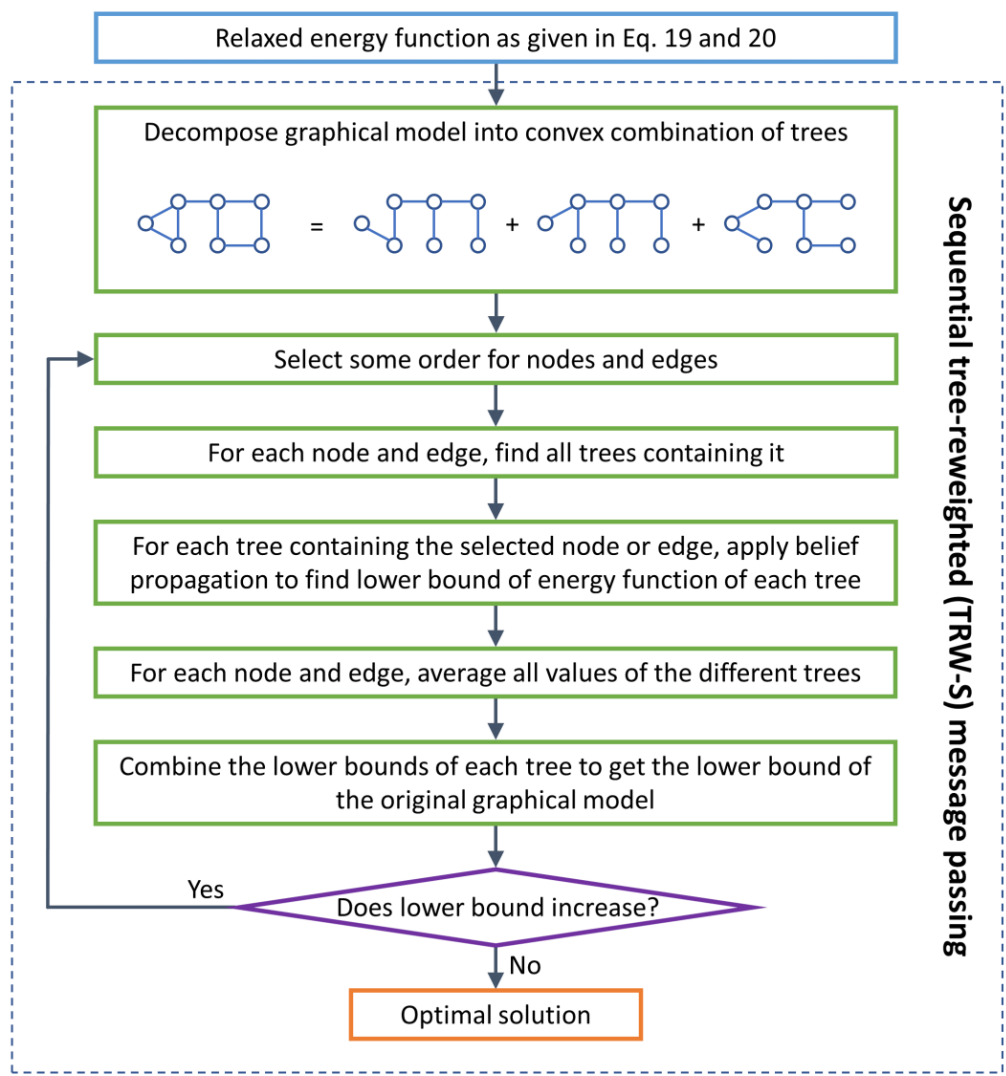

**Figure 10.** Flowchart illustrating the sequential tree-reweighted (TRW-S) message passing algorithm as developed by Kolmogorov [62, 63, 69].

In the next lines, we present the mathematical details of the algorithm. The first step of the TRW-S algorithm consists in decomposing the graphical model into tree-structured functions. The TRW-S algorithm will then try to find a decomposition that maximizes the lower bound of the graphical model (Eq. 19 and 20).

In Sec. 2.2, we showed that the initial wind farm layout optimization problem can be cast as a maximum-a-posteriori problem (Eq. 18), which entails the minimization of the energy function:

$$E(X, \varphi) = \sum_{s \in V} \varphi_s(X_s) + \sum_{(s,t) \in E} \varphi_{st}(X_s, X_t), \quad (21)$$

where $\varphi_s(X_s)$ represents all unary potentials over the model vertices, $V$ (encoding the decision vector), and $\varphi_{st}(X_s, X_t)$ represents the pairwise potentials over the edges, $E$ (encoding the pairwise wake interactions and the turbine number constraint). The minimization problem as given by Eq. 18 can be approximated to the equivalent linear programming (LP) relaxation, Eq. 19 and 20, which TRW-S aims to solve.

We first define the reparameterization: if two parameter vectors (two sets of unary and pairwise potentials) $\varphi$ and $\bar{\varphi}$ define the same energy function, then $\varphi$ is called a reparameterization of $\bar{\varphi}$. We then decompose the graphical model into a convex combination of tree. We consider the unary and pairwise potentials, $\varphi = \{\varphi^T \in \mathbb{R}^d | T \in \mathcal{T}\}$, as a collection of vectors indexed by a finite set $\mathcal{T}$ (set of tree-structured functions) and $\rho$ be a probability distribution on $\mathcal{T}$ (so that $\rho^T \geq 0$ for all $T$ and $\sum_{T \in \mathcal{T}} \rho^T = 1$). The probability distribution is an arbitrary variable that does not affect the optimization results. It is now possible to define a reparameterization of the potentials as follows:

$$\bar{\varphi} \equiv \sum_{T \in \mathcal{T}} \rho^T \varphi^T. \quad (22)$$

The TRW-S algorithm starts with the original parameter vector, $\bar{\varphi}$, and applies iteratively certain reparameterization operations to $\bar{\varphi}$. The goal of these operations is to maximize a lower bound on the energy function, i.e., $\Phi(\bar{\varphi}) = \underset{X \in \{0,1\}^N}{argmin}\, E(X, \bar{\varphi})$ (Eq. 18). Kolmogorov showed that the combined lower bounds of the reparametrized tree-structured functions are equivalent to the lower bound of the energy of the original vector, $\bar{\varphi}$:

$$\sum_{T \in \mathcal{T}} \rho^T \Phi(\varphi^T) \leq \Phi\left(\sum_{T \in \mathcal{T}} \rho^T \varphi^T\right) = \Phi(\bar{\varphi}). \tag{23}$$

The bound described above requires the computation of $\Phi(\varphi^T)$ where the vector $\varphi^T$ corresponds to a tree-structured graph, $T = (V^T, E^T)$. To calculate this bound, the TRW-S algorithm uses max-product belief propagation (BP). The basic operation of BP is passing a message from node $s$ to node $t$ for a directed edge $(s \rightarrow t) \in E^T$. The messages are computed according to:

$$m_{st}(j) = \underset{i \in X_s}{argmin}\{\varphi_s(i) + \varphi_{st}(i,j)\} \quad \forall j \in X_t. \tag{24}$$

The BP algorithm keeps updating the unary and pairwise potentials with the messages according to:

$$\begin{aligned}\varphi_s(j) &:= \varphi_s(j) + m_{st}(j) \\ \varphi_{st}(j) &:= \varphi_{st}(i,j) - m_{st}(j)\end{aligned} \quad \forall j \in X_t, \tag{25}$$

until convergence, which occurs when all edges have valid messages. A valid message indicates that the values $m_{st}(j)$ computed in Eq. 24 satisfy $m_{st}(j) = const_{st}$, where $const_{st}$ does not depend on $j$. The main property of these messages is that they also define a reparameterization.

The TRW-S algorithm works then as follows. First, a reparameterization is conducted on the unary and pairwise potentials of each tree-structured function as the result of the BP algorithm (Eq. 25). This operation finds the minimum energy configuration and minimum energy for each tree. Second, an averaging operation is used to combine the results of each tree according to Eq. 22. The algorithm stops if the value of the lower bound, $\Phi(\bar{\varphi})$ (Eq.23), has not increased within some precision. If a fixed, yet arbitrary, order for the nodes is selected, the TRW-S algorithm guarantees that the lower bound never decreases.